\documentstyle[12pt]{article}
\begin{document}
\author{I. A. Dynnikov and A. P. Veselov}
\title{Integrable Gradient Flows \\ and Morse Theory
\footnote{This work was partially supported by Russian
Foundation for
Fundamental Investigations (grant 94-01-01444), Soros International Science
Foundation (grants MD8000, M3Z000) and INTAS grant (93-0166).}}
\date{}
\maketitle
\begin{center}
{\sf Department of Mathematics and Mechanics\\
Moscow State University\\
119899  Moscow  Russia\\\vspace{2mm}
dynnikov@nw.math.msu.su\hspace{2cm}veselov@nw.math.msu.su}
\end{center}
\bigskip
\begin{abstract}
Examples of Morse functions with integrable gradient flows on some
classical Riemannian manifolds are considered. In particular, we show that
a generic height function on the natural embeddings of classical Lie groups
and certain symmetric spaces is a perfect Morse function, i.~e. has as many
critical points as the homology requires, and the corresponding gradient
flow can be
described explicitly. This gives an explicit cell decomposition and geometric
realization of the homology for such a manifold. As another application of
the integrable Morse functions we give an elementary proof of Vassiljev's
theorem on the flag join of Grassmannians.\end{abstract}
\newcommand{\zamech}[1]{\vspace{2mm}\noindent{\bf Remark #1.}}
\renewcommand{\labelenumi}{\theenumi)}
\newcommand{\ind}{\mathop{\rm ind}\nolimits}
\newcommand{\sgn}{\mathop{\rm sgn}\nolimits}
\newcommand{\eps}{\varepsilon}
\newcommand{\rest}[1]{\Big|_{#1}}
\newcommand{\Frac}[2]{\frac{\displaystyle#1}{\displaystyle#2}}
\newcommand{\eq}[1]{{\rm~(\ref{#1})}}
\newcommand{\ut}[1]{{\rm~\ref{#1}}}
\newcommand{\gr}{\mathop{\rm grad}\nolimits}
\newcommand{\th}{\mathop{\rm th}\nolimits}
\newcommand{\lgr}{{\cal L}}
\newcommand{\dd}{^{\displaystyle\cdot}}
\newcommand{\tr}{\mathop{\rm Tr}\nolimits}
\newcommand{\diag}{\mathop{\rm diag}\nolimits}
\newcommand{\const}{{\rm const}}
\newcommand{\k}{{\bf k}}
\newcommand{\real}{{\bf R}}
\newcommand{\complex}{{\bf C}}
\newcommand{\z}{{\bf Z}}
\newcommand{\quat}{{\bf H}}
\newcommand{\Rep}{\mathop{\rm Re}\nolimits}
\newcommand{\proof}{\noindent{\bf Proof. }}
\newcommand{\proved}{\nolinebreak\mbox{$\triangleright$}}
\newcounter{theorem}[subsubsection]
\newtheorem{theorem}{Theorem}
\newcounter{lemma}[subsubsection]
\newtheorem{lemma}{Lemma}
\newcounter{predl}[subsubsection]
\newtheorem{predl}{Proposition}
\newcounter{sled}[subsubsection]
\newtheorem{sled}{Corollary}
\renewcommand{\thesubsubsection}{\arabic{subsubsection}.}
\renewcommand{\thepredl}{\arabic{subsubsection}.\arabic{predl}}
\renewcommand{\thelemma}{\arabic{subsubsection}.\arabic{lemma}}
\renewcommand{\thetheorem}{\arabic{subsubsection}.\arabic{theorem}}
\renewcommand{\thesled}{\arabic{subsubsection}.\arabic{sled}}
\renewcommand{\theequation}{\arabic{subsubsection}.\arabic{equation}}

\begin{center}{\bf Introduction.}\end{center}
\noindent

Morse theory based on a very natural idea till now remains
one of the most attractive fields in the modern topology. We mention,
for example, the papers by S.~P.~Novikov on a multivalued analogue of
the Morse theory and those by A.T.Fomenko and his colloborators on the
Morse theory for integrable Hamiltonian systems (see e.~g.~\cite{dnf}).

In the present paper we are not concerned with those generalizations but
consider Morse theory in its classical aspects.
It is well-known (see e.~g.~\cite{dnf,%
mil}) that it allows one to extract a topological information
about a manifold $M^n$ from the knowledge of the critical points
of some smooth function $f : M^n \rightarrow {\bf R}$.
If one knows in addition the behaviour of the gradient flow of
the function $f$, it provides a cell decomposition of the manifold $M^n$.

Note that for a
given metric the gradient flow
\begin{equation}
\dot x=\gr f(x)\label{grad}
\end{equation}
is a system of ordinary nonlinear differential equations, which in general
does not admit an explicit integration. The critical points are defined by
the system of nonlinear equations (in algebraic situation - by algebraic one)
\begin{equation}
\label{crit}
\gr f(x)=0,
\end{equation}
which is also difficult to solve in general case.

We will say that  $f$ is an {\it integrable Morse function} on a Riemannian
manifold $M$ if: first, $f$ is a Morse function (or, more generally, Morse-%
Bott function) and, second, the equations\eq{grad} and\eq{crit} are
``integrable''. In the present paper we will mean the integrability in the
naive sense as a possibility of explicit solution leaving apart more serious
discussions of this notion for gradient flows.

We should note that in the literature there exist already examples of gradient
flows which could be naturally considered
as integrable. For example, the Toda
lattice can be interpreted as a gradient flow (see~\cite{Moser, toda}).
Once more remarkable example of an integrable gradient system was recently
pointed out by S.~P.~Novikov in connection with the problem of
``isoperiodical''
deformations for finite-gap potentials and related problems of topological
field theory~\cite{gr-shm, krich}. It has a feature that
the corresponding metric is not positively definite.
For discussions of the integrability of the algebraic equations%
\eq{crit} and examples of such
systems we refer to~\cite{San, A.Ves}.

Such examples naturally appear in the theory of integrable systems.
For instance, the stationary points of the mechanical system with a
potential $U(x)$ are nothing else but
the critical points of $U(x)$: $$\gr U(x)=0.$$ It is natural to expect
that for an integrable system the corresponding stationary problem also has
an integrable nature.

The starting point for us was another example arised from the theory
of integrable systems with discrete time~\cite{v2,mv,v3}.
Such systems are defined by a ``Lagrangian'' which is a function $\lgr$
on the product $Q^{2n}=M^n\times M^n$. A sequence of points from
$M^n$ $X=\{x_k\},\ k\in\bf Z$, is called a solution of the discrete system
with Lagrangian $\lgr$ if $X$ is a critical point of the functional
\begin{equation} S(X)=\sum_{k\in\bf Z}\lgr(x_k,x_{k+1}).\end{equation}

In the papers~\cite{v2,mv,v3} it was shown that in the case when $M^n=SO(N)$
and $\lgr=\tr(XJY^t),\ J=J^t$ is a symmetric matrix, one has an integrable
discrete version of the multidimensional top's dynamics.
Such a form of the Lagrangian is uniquely determined by the following
conditions:

1) symmetry $\lgr(X,Y)=\lgr(Y,X)$;

2) left-invariance $\lgr(gX,gY)=\lgr(X,Y)$, $g\in SO(N)$;

3) bilinearity with respect to $X$ and $Y$ considered in the standard matrix
realization.

Note that because of the left-invariance $\lgr(X,Y)$ depends only on $\omega=
Y^{-1}X:\ \lgr(X,Y)=\lgr(Y^{-1}X,Y^{-1}Y)=\lgr(Y^{-1}X,I)=F(Y^{-1}X)$, where
$F$ is the corresponding function on the group,
\begin{equation} F(\omega)=\tr (J\omega)\label{trj}.\end{equation}
The critical points of the function $F$ correspond to the stationary
solutions of the discrete top's equations.

It turns out that the function\eq{trj} and its analogues on the unitary
and symplectic groups~\cite{v3} can be considered as
examples of integrable Morse functions on the
corresponding group space with biinvariant (Cartan - Killing) metric.
We discuss them in details
in the first and second sections of the present paper. As an application
we get a geometrical representation of the homology classes of these
groups found by V.~A.~Vassiljev in~\cite{vas}.
For example, in the case of the unitary group we have the following
isomorphism:
\begin{equation}\label{vasmakh} H_i(U(n))\cong\bigoplus\limits_{k=0}
^nH_{i-k^2}(G_k(\complex^n)),\end{equation}
where $G_k(\complex^n)$ is the Grassmannian manifold of complex $k$-dimensional
subspaces in $\complex^n$.

Let us define for an arbitrary Schubert cell $\alpha$ in $G_k(\complex^n)$ the
following cycle $[\alpha]$ in $H_*(U(n))$. For each point in $\alpha$, that is
for a $k$-dimensional plane in $\complex^n$, consider the set of all unitary
operators preserving this plane and acting on
its orthogonal complement as the reflection
$(-I)$ (in the paper~\cite{vas} the trivial action on the complement is
considered; this difference is not essential). The union of all such
operators over all points from $\alpha$ forms a cycle $[\alpha]$.
Vassiljev's theorem~\cite[Theorem 1]{vas} states that the set of all cycles
$[\alpha]$ corresponding to all Schubert cells of all Grassmannians
$G_k(\complex^n)$ gives a basis in the integer homology of the unitary
group, which is reflected by equation\eq{vasmakh}. As V.~A.~Vassiljev  told us
this result was found also by M.~Mahowald, so we will call equation%
\eq{vasmakh} {\it Vassiljev-Mahowald decomposition}. We will show
that it can be
obtained easily from the Morse theory for the function\eq{trj}.

We consider also the general height functions
on symmetric embeddings of certain symmetric spaces like
Grassmannians, Lagrangian Grassmannians, the spaces of complex and
quaternian structures and show that they are perfect Morse functions with
integrable gradient flows. For the Grassmannians the corresponding cell
decomposition coincides with the classical Schubert decomposition.

We should note that in a special case $J=\diag(1,0,\dots,0)$ function\eq{trj}
was considered and used for the investigation of the topology of
classical Lie groups
by L.~S.~Pontrjagin in the paper~\cite{p}.

It is interesting that the consideration of another integrable Morse function
allows us to
simplify essentially the proof of the following remarkable result
from the same Vassiljev's paper~\cite{vas}.
Let the Grassmannians
$G_1(V)$, $G_2(V)$,\,$\dots$\,, $G_{n-1}(V)$, where $V$ is an
$n$-dimensional vector space
over $\k=\real,\,\complex$ or
quaternions $\quat$, be embedded into a Euclidean space
$\real^N$ of a large dimension $N$ in a general
position. The union $\Theta_n(\k)$ of
all the $(n-2)$-dimensional simplexes with vertices $(V_1,V_2,\dots,V_{n-1})$,
where
the subspaces $V_i\in G_i(V)$ form a complete flag $F$, that is
$V_1\subset V_2
\subset\dots\subset V_{n-1}$, endowed with the induced topology we will
call the {\it flag join of the Grassmannians}.

The theorem 5 from~\cite{vas} says that this space $\Theta_n(\k)$ is
homeomorphic to the sphere of certain dimension.
In the simplest nontrival case  $\Theta_3(\real)$ can be identified with the
quotient space of ${\bf C}P^2$ by complex conjugation and, therefore,
Vassiljev's result could be considered as higher-dimensional generalization
of the Kuiper-Massey theorem stating that
the last space is homeomorphic to $S^4$. The proof given in~\cite{vas}
is based on some deep results of the modern topology.

In the third paragraph of the present paper we give an elementary proof
of even a little more strong, smooth, version of this theorem. The idea is the
following. Let's consider the set of the symmetric (Hermitian) matrices
satisfying
$$\begin{array}{lll}\tr X&=&0\\\tr X^2&=&1.\end{array}$$
This is a sphere of the dimension we need.
Consider the function $f=\tr X^3$ on it. Its critical points form $(n-1)$
smooth submanifolds isomorphic to $G_1(V)$,\,$\dots$\,, $G_{n-1}(V)$.
The corresponding gradient flow preserves the flag of the eigenspaces
of a matrix
$X$ and determines in a natural way the decomposition of our sphere
into ``flag''
simplexes $\sigma_F$, giving the equivalence with the flag join we need.
 It is interesting to note that the evolution of the eigenvalues
of matrix $X$ is described by the {\it generalized Volterra chain}
(see.~\cite{bog}).

Our paper does not cover the whole variety of the applications of
the integrable Morse functions. Looking at the history of the Morse
theory from this point of view one can find  a number of other remarkable
 examples. For instance, the
well-known proof of the Bott periodicity~\cite{dnf,mil}
uses the functionals
corresponding to integrable systems on the
classical Lie groups.
In the original Bott's proof~\cite{mil} it was
the action functional on the loop
space corresponding to the geodesic flow on the group which is obviously
integrable,
while in Fomenko's version~\cite{Fom} the Dirichlet functional was used,
whose extremals are known in the mathematical physics as the main chiral
fields, the equations of which were integrated by the inverse
scattering method only in 70-th (Zakharov and Mikhailov)~\cite{ts,MZ}.
Of course, in both cases the integrability was used only to find explicitly
some partial solutions, but anyway this also illustrates our approach
to Morse theory.

{}From this point of view it seems to be important that the equations of
the famous {\it self-dual Yang-Mills fields}
can be interpreted as the gradient flow of the Chern-Simons
functional (see e.g.~\cite{t}).
It is well-known that with these equations are connected the recent
remarkable achievements in topology (Donaldson, Witten, Seiberg).

We note also an interesting paper by C.~Tomei~\cite{tomei} in which
the Toda lattice is used for studying topology of manifold of
isospectral Jacoby matrices.

\subsubsection{Height function on a symmetric space as a
perfect Morse function.}
\setcounter{equation}{0}

We start with consideration of the
height functions on the classical Lie groups.
Let $G$ be one of the following compact Lie groups:
$O(n),$ $U(n),$ $Sp(n)$. Here  $Sp(n)$ is the group of quaternian
$n\times n$ matrices $X$ such that $X^{-1}=\overline{X}^T.$ Sometimes it
is called also the unitary symplectic
group because it can be realized as the intersection of the complex
symplectic group $Sp(2n,{\bf C})$ with $U(2n)$.

Consider an arbitrary height function on $G$ embedded in a natural way
into the matrix space $M(n,\k)$, $\k=\real,\,\complex$ or
$\quat$, supplied with the standard scalar product
$$(X,Y) = \Rep\tr(X^*Y).$$ Here for $X\in M(n,\k)$ by $X^*$ we denote
the matrix $\overline{X}^T.$  Such a function has the
form:
\begin{equation}
f_A(X)=\Rep\tr(AX),
\end{equation}
for some matrix $A\in M(n,\k).$

\begin{predl}\label{critpoints}
The critical points of the function $f_A$ are the matrices $X\in G$ 
that satisfy $AX=(AX)^*.$

The gradient flow for the function $f_A$ in the Cartan-Killing metric
on the group $G$ has the form
\begin{equation}\dot X=A^*-XAX.\label{gradspusk}\end{equation}
\end{predl}

\proof Recall that the Cartan-Killing metric on the group $G$ is the
pullback of
the Euclidean structure on the matrix space $M(n,\k)$  under
the embedding. The right hand side of the equation\eq{gradspusk}
is the result of the projection of the vector $A^*$, the gradient of
the function $f_A$ in the ambient
space $M(n,\k)$, onto the tangent space $T_XG$: $X\Big((X^{-1}A^*)-
(X^{-1}A^*)^*\Big)=A^*-XAX$, since $X^{-1}=X^*$ for $X\in G$. The last
expression vanishes iff $AX=(AX)^*$.\proved

\begin{sled}The gradient flow of the function $F_A$ on the group $G$ solves the
problem of the polar decomposition for the  matrix~$A$.
\end{sled}

Indeed, let $$A=JQ,$$ where $J=J^*$ is a positive definite matrix
and $Q\in G$,
is the polar decomposition of the non-degenerate matrix $A$. Then for any
$X$ from $G$ $$\Rep\tr(JX)\le\Rep\tr(J)=\Rep\tr(AQ^*),$$ so, $X=Q^*$ is the
maximum point of the function $f_A$ on $G$. It is not difficult to show
(see also below) that the other critical points are not local maxima.
Therefore, for almost all initial points $X(0)=X_0 \in G$ for the
equation\eq{gradspusk} we will have $X(t)\rightarrow Q^*$
($t\rightarrow\infty$).

\zamech1 It is not necessary to take as a starting point an element of
the group $G$. For example, if we set $X_0=0$ then solving\eq{gradspusk}
we will always have $X(t)\rightarrow Q^*$ ($t\rightarrow\infty$).

To investigate the structure of the set of the critical points of the function
$f_A$, its gradient flow,
and the corresponding cell decomposition of the group $G$ it is
sufficient to consider only diagonal positive matrices
$A=\diag(a_1,\dots,a_n)$ since any nondegenerate matrix can be
reduced to such a form by means of suitable left and right actions
of $G$ corresponding to some rotation of the ambient space $M(n,\k)$.

Let's consider first the case when $A=I$, $I$ is the identity matrix:
$f_I(X)=\Rep\tr X$. It is easy to prove the following result (see for
the details e.g.~\cite{fran}).

\begin{predl}
\label{fI}
\begin{enumerate}
\item
The set of the critical points of the function $f_I$ on $G$ consists of
the involutions: $X^2=I$ and, therefore, can be naturally identified
with the disjoint
union of the Grassmannian manifolds $G_{n,m}(\k),$ $m=0,1,\dots,n$.
Namely, to an arbitrary subspace $V\subset\k^n$ is related the involutive
operator $S_V$, the reflection with respect to $V$:
$$\begin{array}{lll}\vspace{2mm} S_V\rest V&=&\phantom-I\\S_V
\rest{V^\perp}&=&-I,\end{array}$$
and all the involutions can be represented in such a way;
\item
The separatrices of the gradient flow coming from the point $S_V\in G$
fill the domain
with closure isomorphic to $O(m)$, $U(m)$, or $Sp(m)$ correspondingly,
$m=\dim V$, that consists of the operators $X$ such that
$X\rest{V^\perp}=-I$;
\item
the index of a critical point $S_V$, $\dim V=m$, equals to
$$\ind_{S_V}f_A=\left\{\begin{array}{ll}\frac{1}{2}m(m-1)&,\k=\real\\
m^2&,\k=\complex\\m(2m+1)&,\k=\quat\end{array}\right..$$
\end{enumerate}
\end{predl}

Let's  consider now the case when all the eigenvalues of the matrix $A$ are
different.

\begin{theorem} Let $A=\diag(a_1,\dots,a_n),\ 0<a_1<\dots<a_n$. Then
\label{critpnts}
\begin{enumerate}
\item
$f_A$ is a Morse function on the group $G$;
\item
the critical points of the function $f_A$ are the matrices of the form
$$X=\diag(\eps_1,\dots,\eps_n), \mbox{ where } \eps_i=\pm1,\ i=1,\dots,n;$$
\item
the index of such a critical point equals to \begin{equation}\label{ind}
\ind_X(f_A)=\sum_{k=1}^n\delta_{\eps_k,1}((\dim_\real\k) k-1);\end{equation}
\item
the function $f_A$ has the minimal possible number of
critical points among all Morse functions on $G$.
\end{enumerate}\end{theorem}

\noindent{\bf Example.} For
the matrix $$A=\diag(d-1,2d-1,\dots,nd-1),\ \ d=\dim_\real\k$$  the value of
the function $f_A$ at a critical point coincides with the signature
of the Hessian at this point.
This implies that the critical points of the same index lie on the same level
and the greater index the greater level. Such functions are called sometimes
{\it Morse-Smale functions}.

\proof We begin with the item 2). According to proposition\ut{critpoints},
if $X$ is a critical point of $f_A$ then $XAX=A=X^*AX^*=X^{-1}AX^{-1}$,
therefore, $A^2=XA^2X^{-1}$ and the matrices $X$ and $A^2$ commute. It follows
from the diagonal form of the matrix $A$ that such matrix $X$ should be
diagonal too.

Thus, $X_0=\diag (\eps_1,\dots,\eps_n)$, where  $\eps_i=\pm 1$.
Let's compute the second derivative of the function $f_A$ at the point $X_0$.
Consider the exponential map $B\longmapsto X_0\exp(B)$, where $B\in T_IG$.
The diagonal and off-diagonal elements of the matrix $B$ can be taken as
the coordinates in a small neighbourhood of the point $X_0$ under that map.
We have:
$$\begin{array}{l}\vspace{2mm}f_A(X)-f_A(X_0)=\Rep\tr A(X-X_0)=
\Rep\tr AX_0(\exp B-I)=\\\vspace{2mm}\qquad
=\Rep\tr AX_0(B+\frac{1}{2}B^2)+o(\|B\|^2)=\sum\limits_{i,j=1}^na_i\eps_ib_{ij}
b_{ji}+o(\|B\|^2)=\\\qquad
=-\sum\limits_{n\ge j\ge i\ge 1}(a_i\eps_i+a_j\eps_j)|b_{ij}|^2+o(\|B\|^2).
\end{array}$$
Since $a_j>a_i\ge0$ for $j>i$ one has
 $\sgn(a_i\eps_i+a_j\eps_j)=\sgn(a_j\eps_j)$.
Thus $|b_{ij}|^2,\ j\ge i,$ is contained in $d^2f$ with the
negative sign
if $\eps_j=1$ and with the positive one if $\eps_j=-1$. This implies the items
1) and 3).

The item 4) follows from the Morse inequalities~\cite{mil} valid for
an arbitrary Morse function $f$:
\begin{equation}m\label{morse}_j\ge b_j,\end{equation}
where $m_j$ is the number of the critical points of index $j$ of the function
$f$ and $b_j$ is the $j$-th Betti number. For the function $f_A$
a simple calculation shows that the number $m_j$
coincides with the corresponding Betti number
(in the case $G=O(n)$ one should consider the homology groups
over $\z_2$). The Betti numbers of the classical Lie groups have been
calculated first by L.~S.~Pontrjagin in~\cite{p}, one can find them also
e.~g.~in~\cite{fff}.
Thus the Morse inequalities\eq{morse}
for $f_A$ are in fact equalities and this
completes the proof.\proved

A Morse function with the minimal numbers of the critical points is called
a {\it perfect Morse function}.
We will give now a geometrical proof of this property for the height functions,
which does not use the information
about the homology groups of $G$ in a more general situation.

Let  $M$ be a symmetric space embedded in Euclidean Space $E$.
Associate with each point $x\in M$ the reflection $S_x$ of the space $E$
with respect to the normal subspace at the point $x$:
$$S_x(x)=x,\ \ dS_x\rest{T_xM}=-I,\ \ dS_x\rest{(T_xM)^\perp}
=I.$$

We will call an embedding $M\rightarrow E$ {\it symmetric\/} if
all such reflections $S_x$ map $M$ into itself.

The classical groups $G$ with the natural embeddings give the examples
of such symmetric spaces: $M = G$, $E=M(n,\k)$, $S_X(Y)=XY^*X$.
Among other symmetric embeddings\footnote
{All such embeddings are classified, see
\cite{Ferus} and Note added in proof below.} we mention the following:

1. Grassmann manifolds $M=G_{n,m}(\k)$ embedded in the space of symmetric
(Hermitian) matrices $n\times n$ as the set of the reflections with respect
to the corresponding subspaces, $S_X(Y)=XYX$.

2. Lagrangian Grassmannians, consisting of the Lagrangian subspaces in
${\bf R}^{2n} = {\bf C}^n$: $M=U(n)/O(n)$. It is embedded into the space of
complex symmetric matrices as the subset of the unitary matrices,
$X=X^T, X^{-1}=X^*$. The reflection $S_X$ has the
form: $S_X(Y)=X\overline{Y}X$.

3. The symmetric space of complex structures in ${\bf R}^{2n}$:
$M = O(2n)/U(n)$. One can realise it as the set of the skewsymmetric
orthogonal matrices in the Euclidean space of all skewsymmetric matrices,
$S_X(Y) = - XYX$.

4. The space of quaternian structures in ${\bf C}^{2n}$:
 $M = U(2n)/Sp(n)$, realised as the set of the skewsymmetric unitary
matrices in the space of all skewsymmetric matrices,
$S_X(Y) = -X\overline YX$.

5. Symmetric space $M=Sp(n)/U(n)$ can be embedded into the space of
skew-Hermitian quaternian
matrices as the intersection with $Sp(n)$, $S_X(Y)=-XYX$.

The restriction to $M$ of a linear function on $E$ will be referred to as a
{\it height function}.
Among them the Morse functions form an open
everywhere dense set. For example, the function $f_A$ on the classical
Lie group $G$ is a Morse function iff the eigenvalues of the positive symmetric
(Hermitian) matrix $J$ from the polar decomposition $A = JQ$ are pairwise
different (see the calculation above).

\begin{theorem}
\label{SS}
A Morse height function on a symmetric embedding of
symmetric space $M$ is a perfect Morse function.
\end{theorem}

{\bf Proof.}
Notice first of all that by a small variation of such a function $h$
which does not
change the number and indeces of the critical points we can always reduce
to the case when all the critical values of the function $h$ are different and
the gradient separatrix surfaces outcoming from the critical points intersect
the incoming ones transversally. In particular, the separatrices
outgoing from a critical point of index $k$ never hit into critical the
points of indeces $\ge k$ and the set of such separatrices outgoing from a
point of index $k$ and incoming to a point of index $(k-1)$ is finite.

Such a function $h\rest M$
provides a cell decomposition of the manifold $M$ in the following way.
With each critical point of index $k$ we associate a $k$-dimensional cell
which interior is the union of this critical point and all outgoing gradient
separatrices. The incidence number of two cells corresponding to
critical points of indeces
$k$ and $(k-1)$ is defined as the number of gradient separatrices
outgoing from the first point and hitting the second one taken with some sign
depending on the orientation of the cells, which is not essential for us.

Let $x_0$ be one of the critical points of $h\rest M$. Since $dh(T_{x_0}M)=0$
the function $h$ is
invariant under the reflection $S_{x_0}$: $S_{x_0}^*h=h$.
 Thus the set of the critical points
of the function $h\rest M$ is mapped to itself under such reflection.
Since there is no more than one critical point on the same level
of $h$ all the critical points are fixed under $S_{x_0}$.

Consider now the cell complex $C_*$ defined by the function $h\rest M$.
The groups  $C_k$ of this complex consist of the formal sums of the
critical points
of index $k$. Let points $x$ and $y$ have indeces $k$ and $(k-1)$
correspondingly.
As it was shown before, $S_x(y)=y$.
Since $dS_x\rest{T_xM}=-I$\ %
$S_x$ acts freely on the set of the gradient separatrices going from
$x$ to $y$ and from $S_x^2=I$ it follows that the number of these
trajectories is even.

Thus all the incidence numbers of the complex $C_*$ are even and
boundaries of all elements of the complex $C_*\otimes\z_2$ are 0.
According to the Morse theory, the homology of this complex coincides with
$H_*(M,\z_2)$. The fact that the differential of the complex is trivial
implies that Morse inequalities turn out to be equalities and
the number of critical points of the function $h\rest M$ is the minimal
possible one among all Morse functions.\proved

\zamech2 It turns out that for the height functions on
the symmetric spaces the gradient separatrices never hit a critical point
with ``wrong'' index, so we have a correct cell decomposition for {\bf all}
Morse height functions (see the next section).

\zamech3 The height functions turn out to be perfect Morse
functions not only for
symmetric spaces but also for some other classical manifolds.
For instance,
this is the case for the {\it Stiefel manifolds} $V_{n,k}$
naturally embedded into the space of $n\times k$ matrices
and for the {\it flag varieties} realized as the set of symmetric
(Hermitian) matrices with given spectrum.

\subsubsection{Gradient flows and cell decompositions.}
\setcounter{equation}{0}

The functions considered in the previous section have the following remarkable
property: their gradient flows can be linearized by an appropriate change of
variables and, therefore, can be integrated explicitly.

Consider first the gradient flow of the function $f_A(X)=\tr AX,\ A=A^*$
on the group $G=O(n),\ U(n),\ Sp(n)$:
\begin{equation}\dot X=A-XAX.\label{grsp}\end{equation}

Let $$Y=(I-X)(I+X)^{-1}$$ be the {\it Cayley transform} of matrix $X$.
It is known that the Cayley transform is involutive and provides a
bijection between skew-symmetric (skew-hermitian) and
orthogonal (unitary, symplectic) matrices not having $-1$
as an eigenvalue.

\begin{lemma}[\cite{volch}] \label{keli}
The Cayley transform linearizes the flow\eq{grsp}:
\begin{equation}\dot Y=-(AY+YA).\label{ay}\end{equation}\end{lemma}

\proof
$$\begin{array}{l}\vspace{2mm}\dot Y=\Big((I-X)(I+X)^{-1}\Big)\dd=-\dot
X(I+X)^{-1}-(I-X)(I+X)^{-1}\dot X(I+X)^{-1}=\\\vspace{2mm}\qquad
=-(I+Y)\dot X(I+X)^{-1}=-\Frac{1}{2}(I+Y)(A-XAX)(I+Y)=\\\vspace{2mm}\qquad=
-\Frac{1}{2}(I+Y)\Big(A-(I+Y)^{-1}(I-Y)A(I-Y)(I+Y)^{-1}\Big)(I+Y)=\\
\vspace{2mm}\qquad=
-\Frac{1}{2}\Big((I+Y)A(I+Y)-(I-Y)A(I-Y)\Big)=-(AY+YA).\proved\end{array}$$

\begin{predl} The equation\eq{grsp} can be
solved explicitly for arbitrary initial
data $X(0)=X_0$. The solution has the form:\begin{equation}
\label{solution}X(t)=\Big(\sinh(At)+\cosh(At)X_0\Big)\Big(\cosh(At)
+\sinh(At)X_0\Big)^{-1},\end{equation}
where $\sinh$ and $\cosh$ denote the standard
hyperbolic functions of a matrix.
\end{predl}

\proof For an initial matrix $X_0$ not having $-1$ as an eigenvalue the
equation\eq{solution} is obtained by means of the Cayley transform from
the solution of the equation%
\eq{ay}: $Y(t)=\exp(-At)Y_0\exp(-At)$,
with
the initial point $Y_0=(I-X_0)(I+X_0)^{-1}$. The set of the matrices not
having $-1$ in
the spectrum is everywhere dense in $G$,
so the formula\eq{solution}
gives the
solution for all $X_0$ it has a sense.

The function $\th x=\sinh x/\cosh x$ takes values less than 1
for all real $x$, hence the matrix
$\cosh(At)+\sinh(At)X_0=\cosh(At)\Big(I+\th(At)X_0\Big)$ is
non-degenerate for all
$X_0$ from $G$.\proved

\begin{predl}
If the matrices $A_1$ and $A_2$ are symmetric
{\rm(}Hermitian{\rm)} and $A_1A_2=A_2A_1$ then
the gradient flows of the functions $f_{A_1}$
and $f_{A_2}$ on the group $G$
commute.
\end{predl}

\proof
$$\begin{array}{l}\vspace{2mm}
[\gr f_{A_1},\gr f_{A_2}]\Big|_X=\Big(\nabla_{\gr f_{A_1}}(\gr f_{A_2})-
\nabla_{\gr f_{A_2}}(\gr f_{A_1})\Big)\Big|_X=\\\vspace{2mm}\quad=
\nabla_{A_1-XA_1X}(A_2-XA_2X)-
\nabla_{A_2-XA_2X}(A_1-XA_1X)=\\\vspace{2mm}\quad=-(A_1-XA_1X)A_2X-XA_2
(A_1-XA_1X)+\\\vspace{2mm}\qquad+
(A_2-XA_2X)A_1X+XA_1(A_2-XA_2X)=\\\quad=(A_2A_1-A_1A_2)X+X(A_1A_2-A_2
A_1).\end{array}$$
The last expression vanishes for commuting matrices $A_1$ and $A_2$.\proved

\begin{sled}Grassmann manifolds $G_{n,m}(\k)$ embedded into $G$
are invariant under the gradient flow\eq{grsp}. The restricted flow coincides
with the gradient flow of the function $f_A\Big|_{G_{n,m}(\k)}$ with respect
to the induced metric on $G_{n,m}(\k)$.
\end{sled}

\proof Since the gradient flows of functions $f_I$ and $f_A$ commute the last
flow should preserve the set of the first one's stable points. Note that
$G_{n,m}(\k)$ is a connected component of the critical point set of the
function $f_I$ (see.~proposition\ut{fI}). The restricted flow coincides
with the gradient flow of the function $f_A\Big|_{G_{n,m}(\k)}$ with respect
to the induced metric on $G_{n,m}(\k)$ because of the following general
result.\proved

\begin{predl}
\label{restr}
Suppose that a submanifold $N$ of a Riemannian manifold $M$ is invariant
under the gradient flow of some function $f$ on $M$. Then the restricted flow
on $N$ is the gradient one with respect to the restricted function
$f\Big|_{N}$ and induced Riemannian metric on $N$.
\end{predl}

\proof It is easy to see that in general case the gradient flow of the
restricted function $f\rest N$ coinsides with the orthogonal
projection of $\gr f$ in the induced
metric to the tagent subspace of $N$.\proved

It is easy to check that all the symmetric spaces listed in the previous
section are invariant under the gradient flows of the height functions on the
corresponding Lie groups. Indeed, for the Lagrangian Grassmannians
$LG_n = U(n)/O(n)$ this
follows immediately from the flow equation\eq{grsp}.
For the spaces of comlex and quaternian structures $O(2n)/U(n)$
and $U(2n)/Sp(n)$ one can see it from
the equations\eq{gradspusk}:
$$\dot X=A^*-XAX=-(\overline A+XAX)=-\dot X^T$$
since $A^T=-A$ $X^T=-X$.
Similarly, for spaces $Sp(n)/U(n)$ $A^*=-A$, $X^*=-X$,
$$\dot X=A^*-XAX=-(A+XAX)=-\dot X^*.$$

\begin{theorem}
The gradient flows of the height functions on the
classical Lie groups and the embedded symmetric spaces
described above can be integrated explicitly.
\end{theorem}

\proof Using an appropriate shift by an orthogonal (unitary,
symplectic) matrix one can
always reduce the equations of the gradient flow to the form\eq{grsp}
and then apply the explicit formula\eq{solution}.\proved

For the Grassmannians
one can give another explicit formula for the gradient flow.
Let us determine an element $V\in G_{n,m}(\k)$ by $n\times m$  matrix
$Z$ with the columns forming a basis in the subspace $V$. Then the
operator of the reflection with respect to the subspace $V$ has
the form:
$$X=2Z(Z^*Z)^{-1}Z^*-I.$$

\begin{lemma} If $Z(t)$ is changing according to the  equation $$\dot Z=
AZ,$$ then $X(t)$ satisfies\eq{grsp}.\end{lemma}

\proof Indeed, $$\begin{array}{l}\vspace{2mm}\dot X=\Big(2Z(Z^*Z)^{-1}Z^*-I
\Big)\dd=\\\vspace{2mm}\quad=2\dot Z(Z^*Z)^{-1}Z^*+
2Z(Z^*Z)^{-1}\dot Z^*-2Z(Z^*Z)^{-1}(\dot Z^*Z+Z^*\dot Z)(Z^*Z)^{-1}Z^*
=\\\vspace{2mm}\quad=2AZ(Z^*Z)^{-1}Z^*+2Z(Z^*Z)^{-1}Z^*A-
4Z(Z^*Z)^{-1}Z^*AZ(Z^*Z)^{-1}Z^*=\\\quad=A-XAX.\proved\end{array}$$

\begin{predl} The evolution of a subspace $V\in G_{n,m}(\k)$ under the
gradient flow of the function $f_A$ on the Grassmannian $G_{n,m}(\k)$
is given by
 $$V(t)=\exp(At)V.$$
\end{predl}

Let's consider now the cell decompositions determined by our Morse functions.

Let $A=\diag(a_1,\dots,a_n)$, $0<a_1<a_2<\dots<a_n$. As a simple corollary
of the the previous description of the gradient flow one has

\begin{predl} The cell decomposition defined by the function $f_A$ on the
Grassmannian $G_{n,m}(\k)$ coincides with the classical
Schubert decomposition.\end{predl}

The following statement describes the corresponding cell decomposition
of the group $G$, which turns out to be a small modification of the one
considered by V.~Vassiljev~\cite{vas}.

\begin{theorem}\label{schubert}
The cells of the classical Lie group $G$ determined by the Morse function
$f_A$ correspond bijectively
to the Schubert cells of Grassmannians $G_{n,m}(\k)$, $m=0,1,\dots,n$.
The interior of the cell corresponding to a Schubert cell $\sigma$ consists
of the operators such that the eigenspace corresponding to the
eigenvalue $-1$ has the orthogonal complement belonging to
$\sigma$.\end{theorem}

\proof As follows from the proposition\ut{critpnts} all the
critical points of the
function $f_A$ lie on the Grassmannians embedded into $G$. This gives the
bijectivity. It is easy to check that the dimensions of the cells
coincide with the indeces of the corresponding critical points.
Thus
only the invariance of the cells under the gradient flow has to be proven.
Let $X(t)$ be a gradient trajectory. Denote by $V_-(t)$ the eigenspace
of the operator $X(t)$ with the eigenvalue $-1$.
We state that $$V_-(t)=\exp(-At)(V(0)),$$ which implies
$$(V_-(t))^\perp=\exp(At)(V(0))^\perp.$$
Indeed, $$\Big(\cosh(At)+\sinh(At)X(0)\Big)\Big|_{V_-(0)}=
\Big(\cosh(At)-\sinh(At)\Big)\Big|_{V_-(0)}=\exp(-At)\Big|_{(V_-(0))}.$$
This means that $$\Big(\cosh(At)+\sinh(At)X(0)\Big)^{-1}\Big(\exp(-At)
V_-(0)\Big)=V_-(0).$$ Now taking into account\eq{solution} we get:
$$\begin{array}{l}\vspace{2mm}X(t)\Big|_{\exp(-At)V_-(0)}=
\\\vspace{2mm}\qquad=\Big(\sinh(At)+\cosh(At)
X(0)\Big)\Big(\cosh(At)+\sinh(At)
X(0)\Big)^{-1}\Big|_{\exp(-At)V_-(0)}=\\
\qquad=\Big(\sinh (At)-\cosh(At)\Big)\exp(At)
\Big|_{\exp(-At)V_-(0)}=-I\Big|_{\exp(-At)V_-(0)}.\end{array}$$

Thus $\exp(-At)(V(0))\subset V(t)$. In a similar way $V(t)\subset
\exp(-At)(V(0))$. This completes the proof.\proved

We see that the cell decomposition is common for a big family of
the Morse functions $f_A$. This is a particular case of the following
general result.

\begin{theorem}
Any two functions with commuting gradient flows from a connected family of
Morse functions on a compact Riemannian manifold $M$ determine the same cell
decomposition of $M$.
\end{theorem}

\proof From the commutativity it
follows that the gradient flow of the
first function preserves the stationary
set of the second gradient flow. But since they are Morse functions
this is a discrete set and, therefore, both functions have the same critical
points. For given point $x$ the cell containing this point corresponds to
the critical point which is a limit at
the plus
infinity
of the gradient trajectory starting
from $x$ (see above). It is clear that for
functions close enough we have the same limit and, therefore,
the same cell decomposition.
Now the theorem follows from the connectness of the family.\proved

This allows one to know the cell decomposition for such a family
if one knows it for a particular function from it. For example, in the
family $f_A$ we can choose $$A=\diag(d-1,2d-1,\dots,nd-1),\ \ d=\dim_\real\k$$
when the corresponding $f_A$ is a Morse-Smale function (see section 1).

As a corollary of the previous results we have the
Vassiljev-Mahowald decomposition.
\begin{theorem}[\cite{vas}] For the homology of the classical Lie groups one
has the following isomorphisms:
$$\begin{array}{rll}\vspace{2mm}H_i(O(n),\z_2)&\cong&\bigoplus\limits_{k=0
}^nH_{i-\frac{k(k-1)}{2}}(G_{n,k}(\real),\z_2)\\\vspace{2mm}
H_i(U(n),\z)&\cong&\bigoplus\limits_{k=0}^nH_{i-k^2}(G_{n,k}(\complex),\z)\\
H_i(Sp(n),\z)&\cong&\bigoplus\limits_{k=0}^nH_{i-k(2k+1)}(G_{n,k}(\quat),\z).
\end{array}$$\end{theorem}

In a similar way the following result can be proved.

\begin{theorem}
For the homology of Lagrangian Grassmannians, spaces of complex and
quaternian structures and spaces $Sp(n)/U(n)$
one has the following decompositions:
$$
\begin{array}{rll}
\vspace{2mm}
H_i(U(n)/O(n),\z_2)&\cong&\bigoplus\limits_{k=0}^nH_{i-\frac{k(k+1)}2}
(G_{n,k}(\real),\z_2)
\\\vspace{2mm}
H_i(O(2n)/U(n),\z)&\cong&\bigoplus\limits_{k=0}^nH_{i-k(k-1)}
(G_{n,k}({\complex}),\z)
\\\vspace{2mm}
H_i(U(2n)/Sp(n),\z)&\cong&\bigoplus\limits_{k=0}^nH_{i-2k^2}
(G_{n,k}(\quat),\z)
\\\vspace{2mm}
H_i(Sp(n)/U(n),\z)&\cong&\bigoplus\limits_{k=0}^nH_{i-k(k+1)}
(G_{n,k}(\complex),\z).
\end{array}
$$\end{theorem}

For Lagrangian Grassmannians the cell decomposition
determined by a suitable height function turns out to coincide with
one sometimes called Schubert decomposition.\footnote{We would
like to quote D.~B.~Fuks~\cite{top1}: ``The manifold $U(n)/O(n)$ has
a convenient cell decomposition which also could be referred to
as Schubert one although it was invented by Arnold.''}
It could be described as follows (see~%
\cite{top1}): for given Schubert cell $\alpha$ in $G_{n,k}(\real)$
we define the cell $[\alpha] \subset LG_n$ as the set of all
Lagrangian subspaces $L$ in ${\bf C}^n$ with a
standard symplectic structure such that the projection $V$ of $L$ to
${\bf R}^n \subset
{\bf C}^n$ belong to $\alpha$. Indeed, let $U = X + iY$ be a unitary matrix,
which columns represent a basis in $L$. One has $(X + iY)(X^T - iY^T) = I$
and, therefore, $XX^T + YY^T = I$ and $YX^T - XY^T = 0$. This implies that
the symmetric unitary matrix $UU^T$ acts as $-I$ on the orthogonal complement
of $V$ because $UU^T = (XX^T - YY^T) + i(YX^T + XY^T)$ and $X^Tv = 0$ for all
$v \in V^\perp$.

Cell decomposition of the space $Sp(n)/U(n)$ can be constructed in
a way analogous to Latrangian Grassmanians. Geometrically the space
$Sp(n)/U(n)$ can be realized as the set of complex
$n$-dimensional subspaces in $\quat^n$,
on which the restriction of bilinear skewsymmetric function
$\Omega(u,v)=(uj,v)$ is trivial, where
$(\cdot,\cdot)$ is the hermitian scalar product in $\complex^{2n}$.
To any such a subspace $L$ a matrix $U=X+jY$ from $Sp(n)$,
where $X,Y\in M(n,\complex)$, corresponds columns of which
make up a basis in $L$.
Matrix $UiU^*$ is skew-Hermitian (as a quaternian one)
and the operator $-iUiU^*$ acts as $-I$ on $V^\perp$,
where $V\subset\complex^{n}$ is the subspace
spanned over $\complex$ by the columns of matrix $X$.

For the space of complex structures $CS_n$ the gradient cell decomposition
implying isomorphism mentioned above is constructed in the following
way. Let's fix some complex structure $J_0 \in CS_n$, that is an
orthogonal operator in ${\bf R}^{2n}$,
satisfying the relation $J_0^2=-I$. Then for a Schubert cell
$\alpha\subset G_{n,k}(\complex)$ we define $[\alpha]$ as the
union of all complex structures $J$ preserving $V$ and
coinsiding with $J_0$ on $V^\perp$ for all $V \in \alpha$.
In the quaternian case we have a similar construction.

We would like to mention that one can prove other decomposition formulas
using a suitable Morse-Bott height functions with critical submanifolds.
For instanse, using the gradient flow of the height function $f_A$ with
diagonal $A$ having only two different eigenvalues we can prove the
following decomposition formula for the homology of Grassannians.

\begin{theorem}
The following isomorphisms hold for arbitrary $n = n_1 + n_2$:

$$\begin{array}{rll}
\vspace{2mm}
H_i(G_{n,k}(\real),\z_2)&\cong&\bigoplus\limits_{k_1+k_2=k}
H_{i-(n_1-k_1)k_2}(G_{n_1,k_1}(\real)\times G_{n_2,k_2}(\real),\z_2)
\\\vspace{2mm}
H_i(G_{n,k}(\complex),\z)&\cong&\bigoplus\limits_{k_1+k_2=k}
H_{i-2(n_1-k_1)k_2}(G_{n_1,k_1}(\complex)\times G_{n_2,k_2}(\complex),\z)
\\\vspace{2mm}
H_i(G_{n,k}(\quat),\z)&\cong&\bigoplus\limits_{k_1+k_2=k}
H_{i-4(n_1-k_1)k_2}(G_{n_1,k_1}(\quat)\times G_{n_2,k_2}(\quat),\z)
.\end{array}$$
\end{theorem}

A flow similar to ours was used by V.~M.~Buchstaber~\cite[Lemma~3.7]{buch}
in theory of characteristic classes.

\subsubsection{Flag join of Grassmannians and an integrable gradient flow on
the sphere}
\setcounter{equation}{0}

In this section we will give an elementary proof of the following geometrical
result of V.~Vassiljev. Let's consider a general embedding of the
Grassmannians $G_{n,1}(\k)$, $G_{n,2}(\k)$,
$\dots$\,, $G_{n,n-1}(\k)$ into a Euclidean space of a large
dimension. The union of the simplexes with the
vertices $V_1\in G_{n,1}
(\k)$, $V_2\in G_{n,2}(\k)$,\,$\dots$\,, $V_{n-1}\in G_{n,n-1}(\k)$
forming a complete flag: $V_1\subset V_2\subset\dots\subset V_{n-1}$,
with the induced topology we will call the {\it flag join of
Grassmannians\/} and denote following to Vassiljev by $\Theta_n(\k)$.

\begin{theorem}[{V.~Vassiljev~\cite[Theorem 5]{vas}}]
The flag join of Grassmannians is isomorphic to the sphere of the appropriate
dimension:
$$\begin{array}{lll}\vspace{2mm}\Theta_n(\real)&\simeq&S^{\frac{n(n-1)}{2}+n-2}
\\\vspace{2mm}\Theta_n(\complex)&\simeq&S^{n(n-1)+n-2}\\
\Theta_n(\quat)&\simeq&S^{2n(n-1)+n-2}.\end{array}$$
\end{theorem}

To prove this let's consider the function $f=\frac{1}{3}\tr X^3$ on
the set of symmetric
(Hermitian) $n\times n$ matrices $X$ satisfying the conditions:
 \begin{equation}
\label{condition}\begin{array}{lll}\tr X&=&0\\\tr X^2&=&1. \end{array}
\end{equation}
 This is the sphere $S^N$ of dimension
$$ N=n-2+(\dim_\real\k)\frac{n(n-1)}{2}.$$

\begin{predl}\label{grasmtr3}
The set of the critical points of the function $f$ on the sphere $S^N$ is a
disjoint union of the Grassmann manifolds $G_{n,1}(\k),$ $G_{n,2}(\k),
\dots,$ $G_{n,n-1}(\k)$ smoothly embedded into $S^N$ in such a way that to
a subspace $V\in
G_{n,m}(\k)$ corresponds the matrix $X=X(V)$ such that \begin{equation}
\begin{array}{llll}\vspace{2mm}X\rest V&=&-\sqrt{\frac{n-m}{nm}}&I
\\X\rest{V^\perp}&=&\phantom{-}\sqrt{\frac{m}{n(n-m)}}&I
\label{V->X}.
\end{array}\end{equation}\end{predl}

\proof Functions $\tr X$, $\tr X^2$, $\tr X^3$ in the space of symmetric
(Hermitian) matrices have gradients: $$\begin{array}{lll}\vspace{2mm}
\gr(\tr X)&=&I\\
\vspace{2mm}\gr(\tr X^2)&=&2X\\\gr(\tr X^3)&=&3X^2.\end{array}$$

The critical points of the function $f$ are those $X$ for which matrices
$I$, $X$, $X^2$ are dependent. This holds iff $X$ has no more than two
different eigenvalues. Since $\tr X=0$, $\tr X^2=1$ such a matrix $X$
has exactly
two different eigenvalues $\mu_1$ and $\mu_2$, $\mu_1<0<\mu_2$. Let $m$ and
$(n-m)$ be their multiplicity. It is easy to find from\eq{condition}
that $$\mu_1=-\sqrt{\frac{n-m}{nm}},\ \ \ \mu_2=\sqrt{\frac{m}{n(n-m)}}.$$
An operator with such a spectrum is uniquely determined by its eigenspace
corresponding to $\mu_1$.\proved

\begin{lemma} The gradient flow of the function $f$ has the form:
\begin{equation}\label{grsptr3}\dot X=X^2-f(X)X-\frac{1}{n}I.\end{equation}
\end{lemma}

\proof The right hand side of the equation\eq{grsptr3} is the orthogonal
 projection of
the vector $X^2$ to the tangent
space of our $S^N$ with the normal vectors $X$
and $I$: $$X^2-\frac{(X^2,X)}{(X,X)}X-\frac{(X^2,I)}{(I,I)}
I=X^2-\frac{\tr X^3}{\tr X^2}X-\frac{\tr X^2}{\tr I}I=X^2-(\tr X^3)X-
\frac{1}{n}I.$$ It is well-known that it gives the gradient vector of
the restricted
function with respect to the induced metric.\proved

\begin{predl} The gradient flow of the function $f$ preserves the eigenspaces
and the anharmonic ratios
of the eigenvalues of the matrix $X$.\label{sl2}\end{predl}

\proof The first statement follows from the fact that the left hand side of the
equality\eq{grsptr3} is a polynomial of $X$ with scalar coefficients.
The dynamics of the eigenvalues is determined by a vector field of
the form $(a \lambda^2+b \lambda+c)\frac{\partial}
{\partial \lambda}$, where $(a,b,c)\ne(0,0,0)$, which is the same for all of
them. But it is well-known that such a vector field generates a one-parametric
group of the projective transformations of the line.\proved

Let \begin{equation}\label{vozr}\lambda_1(X)\le\lambda_2(X)\le\dots\le
\lambda_n(X)\end{equation} be the eigenvalues of the matrix
$X$, $$V_1(X),\ V_2(X),\dots,\ V_n(X)$$ be the corresponding eigenspaces.
Let's define the {\it eigenflag} of $X$ as
\begin{equation}\label{flag}F(X)=\{U_0\subset U_1\subset\dots\subset
U_n\},\end{equation} where
$$U_0(X)=\{0\},\ \ U_k(X)=\bigoplus\limits_{i=1}^kV_i(X),\ \ 1\le k\le n.$$
If some eigenvalues are equal to each other
the corresponding eigenflag will be incomplete.

Proposition\ut{sl2} gives the set of integrals of the
flow\eq{grsptr3}, consisting of the eigenflag of $X$ and the anharmonic ratios
of the eigenvalues.
This set is complete in the sense that if we fix the values of the integrals
we will have a trajectory.

Indeed, a symmetric (Hermitian) matrix is completely determined by its
spectrum and
the eigenflag. Given the anharmonic ratios
the spectrum is defined up to a
projective transformation. The group of projective transforms of the line is
3-dimensional, therefore, the conditions
\begin{equation}\label{lambdas}\begin{array}{lll}
\sum\limits_{i=1}^n\lambda_i&=&0\\\sum\limits_{i=1}^n\lambda_i^2
&=&1,\end{array}\end{equation} equivalent to\eq{condition} determine a
1-dimensional family.

For two flags $F_1$ and $F_2$ we will write $F_1\le F_2$
if all the subspaces of $F_1$ are contained also in $F_2$.

\begin{lemma}For a given flag $F=\{U_0=\{0\}\subset U_1\subset U_2\subset
\dots\subset U_{k+1}=\k^n\}$ the matrices $X$ satisfying\eq{condition}
and such that $F(X)\le F$ form a
$(k-1)$-dimensional simplex with the vertices corresponding
to the subspaces $U_1,\dots,U_k$  according to\eq{V->X}. The ratios
\begin{equation}a_i=\frac{\lambda_{i+1}
(X)-\lambda_i(X)} {\lambda_n(X)-\lambda_1(X)},\label{a}
\end{equation}are the barycentric coordinates on this simplex.
\end{lemma}

\proof For simplicity we give the proof only for the case of complete flag $F$.
In this case we can assume without loss of generality that the subspaces
$U_1,\dots,U_n$ are ${<}e_1{>},{<}e_1,e_2{>},\,%
\dots\,,{<}e_1,\dots,e_n{>}$, where $e_1,\dots,e_n$ is the standard basis in
$\k^n$. Then the condition $F(X)\le F$ is equivalent to the statement that $X$
is a diagonal real matrix with the increasing elements at the diagonal:
$X_{11}\le X_{22}\le\dots\le X_{nn}$.

So, there is one-to-one correspondence between the set of matrices $X\in S^N$
satisfying $F(X)\le F$ and the set of sequences
$\lambda_1\le\dots\le\lambda_n$, satisfying\eq{condition}.
We should only check that the functions\eq{a} can be considered
as the barycentric coordinates on this set.

Indeed, by construction we have $$\begin{array}{c}\vspace{2mm}\displaystyle
\sum\limits_{i=1}^{n-1}
a_i=1,\\a_i\ge0\mbox{   for all }i=1,\dots,n-1.\end{array}$$
The sequence $\{\lambda_i\}$ is defined by $\{a_i\}$
from\eq{a} uniquely up to an affine
transform. The normalization conditions\eq{lambdas} together with the
increasing property\eq{vozr} determine it uniquely.\proved

Summarizing, we have the following \begin{theorem}
The gradient flow of the function $f = \frac{1}{3}tr X^3$ on the sphere
in the space of symmetric matrices $X$: $tr X = 0,\ tr X^2 = 1$, is integrable
and gives a decomposition of this sphere into a union of simplexes. The
vertices of such a  simplex belong to the Grassmannians embedded into
the sphere as the critical set of $f$ and form a flag.
\end{theorem}

Concerning the integrability we would like to mention
that the gradient flow of the function $f$ in the
barycentric variables $a_i$ after a suitable change of time has the
form
of the {\it generalised Volterra chain} (see~\cite{bog}):
\begin{equation}\label{volter}\frac{d}{d\tau}a_i=a_i
\left(
\sum\limits_{k=i-n}^{i-1}a_k-\sum\limits_{l=i+1}^{i+n}a_l
\right),
\end{equation}
where we set $a_i$ with $i\le 0$ and $i>n$ to be zero.

Indeed, from\eq{grsptr3} we have: $$\dot\lambda_i=
\lambda_i^2-I_3\lambda_i-\frac {1}{n},\mbox{ where
}I_3=\sum\limits_{k=1}^n\lambda^3.$$

Hence
$$\begin{array}{l}\vspace{4mm}\dot a_i=\left(
\Frac {\lambda_{i+1}-\lambda_i} {\lambda_n-\lambda_1}
\right)\dd=
\Frac {\dot\lambda_{i+1}-\dot\lambda_i} {\lambda_n-\lambda_1}-
\Frac {(\dot\lambda_n-\dot\lambda_1)(\lambda_{i+1}-\lambda_i)}
{(\lambda_n-\lambda_1)^2}=
a_i\left(
\Frac {\dot\lambda_{i+1}-\dot\lambda_i} {\lambda_{i+1}-\lambda_i}
-\Frac {\dot\lambda_n-\dot\lambda_1} {\lambda_n-\lambda_1}
\right)=\\\vspace{4mm}\phantom{a_i}=
a_i\left(
\Frac {\lambda_{i+1}^2-I_3\lambda_{i+1}-(\lambda_i^2-I_3\lambda_i)}
{\lambda_{i+1}-\lambda_i}-
\Frac {\lambda_n^2-I_3\lambda_n-(\lambda_1^2-I_3\lambda_1)}
{\lambda_n-\lambda_1}
\right)=\\\phantom{a_i}=
a_i(\lambda_{i+1}+\lambda_i-\lambda_n-\lambda_1)=
(\lambda_n-\lambda_1)a_i
\left(
\displaystyle\sum\limits_{k=1}^{i-1}a_k-\sum\limits_{l=i+1}^{n-1}a_l
\right).
\end{array}$$
Let's change now the time variable on the gradient flow trajectory:
$d\tau=(\lambda_n
-\lambda_1)\,dt$. The equation will take the form:
\begin{equation}\frac{d}{d\tau}a_i=a_i
\left(
\sum\limits_{k=1}^{i-1}a_k-\sum\limits_{l=i+1}^{n-1}a_l
\right),\label{volt}
\end{equation}
which coincides with\eq{volter} if one set $a_0=a_{-1}=a_{-2}=\dots=a_n
=a_{n+1}=a_{n+2}=\dots=0$.

In this case the generalized Volterra chain can be easily explicitly
integrated. To do this let us change variables (cf.~\cite{bog}):
$$b_i=\sum\limits_{k=1}^ia_i=\frac{\lambda_{i+1}-
\lambda_1}{\lambda_n-\lambda_1}.$$ We will have:$$\frac{d}{d\tau}b_i=
b_i\frac{\lambda_{i+1}+\lambda_1-\lambda_n-\lambda_1}{\lambda_n-\lambda_1}=
b_i(b_i-1),$$
hence $$b_i=\frac{1}{1-c_ie^t},$$
where the constants $c_i$ are determined by the initial data:
$c_i = 1 - b_i(0)^{-1}$.
The last formula is in a good agreement with our previous considerations.

As a corollary we have the following strengthened version of Vassiljev's
theorem.

\begin{theorem} Flag join of the Grassmannians $\Theta_n(\k)$ admits
the smooth structure of the standard sphere for which all the
embeddings $G_{n,m}\hookrightarrow \Theta_n(\k)$ are smooth.\end{theorem}

{\bf Acknowledgements.}
We would like to thank V.~M.~Buchstaber for pointing out the references%
{}~\cite{fran,buch}.
One of us (A.~V.) is grateful to J.~Moser and P.~Santini for
useful discussions and to Forschungsinstitut f\"ur
Mathematik (ETH,~Zurich) for the
hospitality during the summer term of 1995.

\vspace{1cm}
{\bf Note added in proof.} Recently L.~Nicolaescu attracted our attention to
the papers \cite{Take1,Take2,Ram}, where the properties of
the height Morse functions on a class of so-called {\it R-symmetric spaces}
 are investigated. These spaces turn out to be precisely the ones for which
embeddings with the property we demanded in section 1 do exist \cite{Ferus}.
It is interesting that the same embeddings are remarkable from
another geometrical
point of view: they are minimal in the sense of total curvature (see~%
\cite{Koba,TaKo}). The whole list of R-symmetric spaces besides the examples we
mentioned in section 1 consists of $M = Sp(n)/U(n)$ and exceptional
spaces $M = E_6/Spin(10) \times T, F_4/Spin(9)$ or Cayley projective plane,
$E_7/E_6 \times T$ and $E_6/F_4$ (see \cite{Koba}). It would be interesting
to investigate the integrability of the gradient flows for the height
functions on the corresponding embeddings.
We are very grateful to L.Nicolaescu and E.Leuzinger for pointing out these
remarkable papers to us.

\end{document}